 \documentclass[11pt]{article}
 \usepackage[top=1in,bottom=1in,left=1.5in,right=1.5in]{geometry}
 \usepackage{amsmath,amssymb}
  \usepackage{latexsym}
  \usepackage{graphicx}
  \usepackage{epsfig}
  \usepackage{subfigure}

 \def\cR{{\cal R}}

 \def\bA{{\bf A}}

 \def\c{{\bf c}}
 
 \def\e{{\bf e}}

 \def\1{{\bf 1}}

 \def\o{{\bf o}}

 \def\bt{{\bf t}}
 \def\u{{\bf u}}
 \def\v{{\bf v}}
 \def\w{{\bf w}}
 \def\x{{\bf x}}
 \def\y{{\bf y}}

 \def\X{{\bf X}}
 \def\y{{\bf y}}
 
 \def\diag{\mathop{{\rm diag}}\nolimits}
 \def\hs{\hspace*{\parindent}}

 \newcommand{\qed}{\hspace*{\fill} $\Box$\\}
 \def\rF{{\rm F}}
 
 \def\rO{{\rm O}}
 \def\rS{{\rm S}}
 \def\rank{\mathop{{\rm rank}}\nolimits}
 \def\span{\mathop{{\rm span}}\nolimits}
 
 \def\t{^{\rm T}}
 \def\tr{\mathop{{\rm trace\;}}\nolimits}

 \def\R{\mathord{\mathbb R}}

 \newtheorem{theo}{\bf \hs Theorem}[section]

 \newtheorem{corol}[theo]{\bf \hs Corollary}

\begin{document}

 \title{Fast Monte-Carlo Low Rank Approximations for Matrices}
 \author{Shmuel Friedland, Amir Niknejad
 \\Department of Mathematics, \\ Statistics and Computer
 Science\\ University of Illinois at Chicago
 \\ Chicago, Illinois 60607-7045\\
 friedlan@uic.edu, niknejad@math.uic.edu
 \and Mostafa Kaveh, Hossein
 Zare  \\Department of Electrical\\and Computer Engineering \\
 University of Minnesota\\Minneapolis, MN 55455\\
 \{mos, hossein\}@ece.umn.edu}

 \date {  }


 \maketitle

 \begin{abstract}
 In many applications, it is of interest to approximate data, given
 by $m\times n$ matrix $A$, by a matrix  $B$ of at
 most rank $k$, which is much smaller than $m$ and $n$.
 The best approximation is given by singular value decomposition,
 which is too time consuming for very large $m$ and $n$.

 We present here a Monte Carlo algorithm for iteratively computing a $k$-rank
 approximation to the data consisting of $m\times n$ matrix $A$.
 Each iteration involves the reading of $O(k)$ of columns
 or rows of $A$.  The complexity of our algorithm is
 $O(kmn)$.
 Our algorithm, distinguished from other known algorithms,
 guarantees that each iteration is a better $k$-rank
 approximation than the previous iteration.  We believe that this
 algorithm will have many applications in data mining, data storage and
 data analysis.

 \noindent

 \end{abstract}


 \section{Introduction}
 In many applied settings, dealing with a very large data set, it is
 important to reduce the size of data in order to
 make an inference about the features of data set, in a timely
 manner. Assume that the data set is represented by an $m\times n$ matrix
 $A\in \R^{m\times n}$.  It is important to find an approximation
 $B\in\R^{m\times n}$ of a
 specified rank $k$ to $A$, where $k$ is much smaller than
 $m$ and $n$.

 Here are several motivation to obtain such $B$.
 First, the
 storage space needed for $B$ is $k(m+n)$, which is much smaller
 than the storage space $mn$ needed for $A$.
 Indeed, $B$ can be represented as
 \begin{equation}\label{bexpres}
 B=\x_1\y_1\t+\x_2\y_2\t+\ldots\x_k\y_k\t,
 \end{equation}
 where $\x_1,\ldots,\x_k$ and $\y_1,\ldots,\y_k$ are $k$ column vectors with
 $m$ and $n$
 coordinates respectively.  We store the vectors
 $\x_1,\ldots,\x_k,\y_1,\ldots,\y_k$, which need the storage
 $k(m+n)$, and compute the entries of $B$, when needed, using the
 above expression of $B$.

 The second most
 common application is clustering algorithms
  as in \cite{APWYP}, \cite{AGGR},
 \cite{DFKVV}, \cite{EKSX} and \cite{PAJM}.
 Assume that our data represents $n$ points and each point has $m$
 coordinates.  That is each point represented by a column of the
 matrix $A$.  We want to cluster the points in such a way that the
 distance between points in the same cluster is much smaller than the
 distance between any two points from different clusters.
 One way to do this is to project all the point on the $k$ main
 orthonormal directions encoded by the first $k$ left singular
 vectors of $A$.  Then cluster using either $k$ one dimensional
 subspaces or the whole $k$ dimensional subspace.  The $k$-rank
 approximation $B$ gives the approximation to this $k$ dimensional
 subspace and the approximation to the first $k$ singular vectors.
 Another way to cluster is to use the projective clustering.
 It is known that a fast SVD, i.e. fast $k$-rank approximation,
 is a main tool in this area \cite{AGGR} and \cite{DFKVV}.

 The third application is in DNA microarrays, in particular
 in data imputation \cite{FNC}.  Let $A$ be the gene expression
 matrix.  The $m$ rows of the matrix $A$ are indexed by the genes,
 while the $n$ columns are indexed by the number of experiments.
 It is known that the effective rank of $A$ is $k$, is
 usually much less then $n$.  The SVD decomposition is the main
 ingredient of the FRAA, (\emph{fixed rank approximation algorithm}),
 which successfully implemented in \cite{FNC} to impute the corrupted
 entries of $A$.  The fast $k$- approximation algorithm suggested
 in this paper, combined with the clustering of similar genes, can be
 implemented to improve the FRAA algorithm.

 The best approximation
 $B$, which minimizes the Frobenius norm $||A -B||_{\rF}^2:=\tr
 (A-B)\t (A - B)$, is given by the
 \emph{singular value decomposition} (SVD) of $A$.
 However SVD decomposition needs $O(mn\cdot \min(m,n))$ time computation,
 which is often too prohibitive.

 One way to find a \emph{fast} $k$-rank approximation is to
 choose at random $l\ge k$ columns or rows of $A$
 and obtain from them $k$-rank approximations of
 $A$.  This is basically the spirit of the algorithm suggested
 in \cite{FKV}.  We call this algorithm the FKV algorithm.
 Assuming a statistical model for the
 distribution of the entries of $A$ the authors give some error
 bounds on their $k$-rank approximation.  The weak point of FKV
 algorithm is its inability to improve iteratively FKV
 approximation by incorporating additional parts of $A$.
 In fact in the recent paper \cite{DFKVV}, which uses FKV random
 algorithm, this point is mentioned specifically in the end of the paper:
 ``..., it would be
 interesting to design an algorithm that improves this
 approximation by accessing $A$ (or parts of $A$) again."

  The aim of this paper is to provide a sampling framework for
 iterative updates of k-rank approximations of $A$, by reading
 iteratively additional
 columns or rows of $A$, which improves \emph{for sure} the
 approximation $B$ each time it is updated.  The quality of the
 approximation of $B$ is given by the Frobenius norm $||B||_{\rF}$,
 which is a nondecreasing sequence under these iterations.
 The rate of increase of the norms $||B||_{\rF}$ can be used as a
 stopping rule for terminating the algorithm.  Also the updating
 algorithm of $k$-rank approximation gives approximation to the
 first $k$ singular values of $A$, and the approximations to the
 first $k$ left and right singular vectors of $A$.

 Assuming that the number of columns or rows which are read is
 $l=O(k)$, the complexity of our algorithm is $O(kmn)$.
 The intensive part of the computations is devoted to obtain
 the spectral decompositions of $(k+l)\times (k+l)$  real symmetric
 matrices, where the computation for each decomposition is of order $O(k^3)$.
 The $O(kmn)$ part of the algorithm is due to the multiplications
 of $k$ column vectors, which approximate the $k$ left or right
 singular eigenvectors of $A$, by $A\t$ or $A$ respectively.
 This part of the algorithm can be parallelized, which
 will speed significantly the algorithm suggested here.
 The simulations that we performed show that we need a small
 number of iterations, (around 5), to get a very good
 approximation to the best $k$-rank approximation of $A$.

 We believe that this
 algorithm will have many applications in data mining, data storage and
 data analysis.

 \section{SVD}

 \setcounter{equation}{0}

 In this section, we recall some basic facts about SVD,
 which are embedded in our algorithm, see \cite{GV}. Let
 $\R^m,\R^{m\times n},\rO_{md}(\R),
 \rS_n(\R),
 $ be the linear space of real
 column vectors with $m$ coordinates, the
 linear space of real $m\times
 n$ matrices, the subset of $m\times d$ real valued matrices whose $d \;(\le m)$
 columms is an orthonormal system and
 the subspace of $n\times n$ real
 symmetric matrices.
 For $S\in \rS_n(\R)$ we
 let $S\ge 0$ if $S$ is nonnegative definite and $S>0$ if $S$
 is positive definite.
 Denote by $\diag(d_1,\ldots,d_m)\in\R^{m\times m}$
 the matrix with
 the diagonal entry $d_i$ on the $(i,i)$ position for
 $i=1,...,m$ and all other entries are equal to zero.
 Let $r=\rank A$ be the rank of $A$.  Then the SVD
 decomposition of $A$ is given as
 $A=U_r \Sigma_r V_r\t$ where $U_r\t U_r=V_r\t V_r=I_r$,
 and $\Sigma_r=\diag(\sigma_1,...,\sigma_r),\;
 \sigma_1\ge...\ge\sigma_r > 0$.
 Let $\u_1, \cdots, \u_r\in\R^m$ and  $\v_1,\dots,\v_r\in\R^n$
 denote $1,\ldots,r$ orthonormal columns of $U$ and $V$ respectively.
 Then  $A=U_r \Sigma_r V_r\t$ can be written as
 $A=\sum_{q=1}^r \sigma_q\u_q\v_q\t$.  The vectors $\u_i,\v_i$
 are called the
 \emph{left} and \emph{right singular vectors} of $A$
 respectively, which
 correspond to the singular value $\sigma_i$.
 The left and the right singular vectors can be computed from the
 right and the left singular vectors by the formulas:
 \begin{equation}\label{rleftsvec}
 \u_i=\frac{1}{\sigma_i}A\v_i, \quad
 \v_i=\frac{1}{\sigma_i}A\t\u_i, \quad i=1,\ldots,r.
 \end{equation}
 Equivalently, $\sigma_1^2\ge\ldots\ge\sigma_r^2>0$ are all the
 positive eigenvalues of the nonnegative definite symmetric
 matrices $A A\t,A\t A$, with the corresponding orthonormal eigenvectors
 $\u_1, \cdots, \u_r\in\R^m$ and  $\v_1,\dots,\v_r\in\R^n$.

 Denote by $||A||_{\rF}:=\sqrt{\tr A\t A}=\sqrt{\tr A A\t}$
 the Frobenius ($\ell_2$) norm of $A$. It
 is the Euclidean norm of $A$ viewed as a vector with $nm$
 coordinates. Each term $\u_q\v_q\t$ in SVD decomposition of $A$ is a rank one
 matrix with $||\u_q\v_q\t||_{\rF}=1$. Let $\cR(n,m,k)$ denote the
 set of $n\times m$ matrices of at most rank $k$ ($\min(m,n)\ge k$).  Then
 for each $k$,  $k\le r$, the SVD of $A$ gives the solution to the
 following approximation problem:

 \begin{equation}\label{mincharsvd}
 \min_{B\in\cR(n,m,k)} ||A-B||_{\rF}=||A-\sum_{q=1}^k
 \sigma_q\u_q\v_q\t||_{\rF}= \sqrt{\sum_{q=k+1}^r\sigma_q^2}.
 \end{equation}

 If $\sigma_k > \sigma_{k+1}$  then $\sum_{q=1}^k
 \sigma_q\u_q\v_q\t$ is the unique solution to the above minima
 problem. For the purposes of this paper, it will be convenient to
 assume that $\sigma_q=0$ for any $q >\rank A$.

 Given $k\in [1,r]$ then $\sigma_1^2,\ldots,\sigma_k^2$,
 $\u_1,\ldots,\u_k$ and $\v_1,\ldots,\v_k$ are characterized by the
 following maximal characterization:

 \begin{theo}\label{maxchar}  Let $A\in \R^{m\times n}$ and $k\in
 [1,\min(m,n)]$.  Let $\x_1,\ldots,\x_k$ and
 $\y_1,\ldots,\y_k$ be two sets of orthonormal vectors
 in $\R^m$ and $\R^n$ respectively.  Then
 \begin{equation}\label{maxchar1}
 \sum_{i=1}^k (A\t\x_i)\t (A\t\x_i)\le \sum_{i=1}^k
 \sigma_i^2,\;
 \sum_{i=1}^k (A\y_i)\t (A\y_i)\le \sum_{i=1}^k \sigma_i^2.
 \end{equation}
 Equality in the first or second inequality occurs if
 $\span(\x_1,\ldots,\x_k)$ or $\span(\y_1,\ldots,\y_k)$, contains $k$
 linearly independent left or right singular vectors of $A$,
 corresponding the the $k$ maximal singular values of $A$ respectively.

 \end{theo}

 The above characterization follows from the maximal, (Ky Fan
 characterization), of the sum of the first biggest eigenvalues of a
 real symmetric matrix:
 \begin{theo}\label{kyfanchar}  Let $S\in S_p(\R)$ be a real
 $p\times p$ symmetric matrix.  Let $\lambda_1\ge\ldots\ge\lambda_p$
 be the eigenvalues $S$ arranged in a decreasing order and listed
 with their multiplicities.  Let $\w_1,\ldots,\w_p\in\R^p$ be an
 orthonormal set of the corresponding eigenvectors of $S$:
 $S\w_i=\lambda_i\w_i, i=1,\ldots,p$.  Let $k\in [1,p]$ be an integer.
 Then for any orthonormal set $\x_1,\ldots,\x_k\in\R^p$
 \begin{equation}\label{kyfanin}
 \sum_{i=1}^k \x_i\t S\x_i \le \sum_{i=1}\lambda_i=\sum_{i=1}^k \w_i\t S\w_i.
 \end{equation}
 Equality holds if and only if the subspace
 $\span(\x_1,\ldots,\x_k)$ contains $k$ linearly independent
 eigenvectors of $S$ corresponding to the eigenvalues
 $\lambda_1,\ldots,\lambda_k$.
 \end{theo}

 See for example \cite{Fr77} for proofs and the
 references.
 Note that $\x_1,\ldots,\x_k\in\R^m$
 is a system of orthonormal vectors in
 $\R^m$ if and only if the $m\times k$ matrix $(\x_1,\ldots,\x_k)$
 is in $\rO_{mk}(\R)$.

 To obtain Theorem \ref{maxchar} from Theorem \ref{kyfanchar}
 we let $p=m$, (or $p=n$), and $S$ to be equal to $A A\t$, (or $A\t
 A$).
 In (\ref{maxchar1}) we emphasized the
 complexity of the computations of the left-hand side of the
 inequalities.  See also \cite{FNC} for applications of Theorem
 \ref{kyfanchar} to data imputation in DNA microarrays.

 In the rest of the paper we give the version of our results
 to $k$-orthornormal
 systems $\x_1,\ldots,\x_k\in\R^m$.  Similar results holds
 for $k$-orthonormal systems $\y_1,\ldots,\y_k\in\R^n$.

 \begin{corol}\label{minchar}  Let $A\in\R^{m\times n}$ and $k\in
 [1,\min(m,n)]$ be an integer.  Then for any $k$-orthornormal
 system $\x_1,\ldots,\x_k\in\R^m$
 the following equality holds:
 \begin{equation}
 ||A-\sum_{i=1}^k \x_i (\x_i\t A)||^2_{\rF}=||A||^2_{\rF} -
 \sum_{i=1}^k (A\t\x_i)\t (A\t\x_i). \label{id1}
 \end{equation}
 In particular the best $k$-rank approximation of $A$ is given by
 $\sum_{i=1}^k \u_i (\u_i\t A)$,
 where $\u_1,\ldots,\u_k\in \R^m$ is an orthonormal sets of the left
 singular vectors of $A$ corresponding to
 $\sigma_1,\ldots,\sigma_k$.
 \end{corol}

 The next theorem is the key theorem for updating the $k$-rank
 approximation for

 \noindent
 $\sum_{i=1}^k \x_i (\x_i\t A)$, for some $(\x_1,\ldots,\x_k)\in
 \rO_{mk}(\R)$.

 \begin{theo}\label{updatel}  Let $\x_1,\ldots,\x_k\in\R^m$,
 be an orthonormal system in $\R^m$.
 Let $\w_1,\ldots,\w_l\in\R^m$, be a given set in $\R^m$.
 Perform the modified Gram-Schmidt
 algorithm on $\x_1,\ldots,\x_k,\w_1,\ldots,\w_l$,
 to obtain an orthonormal
 set $\x_1,\ldots,\x_p\in\R^m$,
 where $k\le p\le k+l$.  Assume that $k < p$, i.e. $\span
 (\w_1,\ldots,\w_l)\nsubseteq \span (\x_1,\ldots,\x_k)$.
 Form $p\times p$ real symmetric matrix $S:=((A\t\x_i)\t
 (A\t\x_j))_{i,j=1}^p$, and assume that
 $\lambda_1\ge\ldots\ge\lambda_k$ are the $k$-largest eigenvalues of
 $S$ with the corresponding $k$-orthonormal vectors
 $\o_1,\ldots,\o_k\in \R^p$.  Let $O:=(\o_1,\ldots,\o_k)\in
 \rO_{pk}(\R)$ and define $k$-orthonomal vectors
 $\tilde \x_1,\ldots,\tilde \x_k\in \R^m$ as follows:
 \begin{equation}\label{xtildef}
 (\tilde \x_1,\ldots,\tilde \x_k)=(\x_1,\ldots,
 \x_p)O.
 \end{equation}
 Then
 \begin{equation}
 \sum_{i=1}^k (A\t\x_i)\t (A\t\x_i)\le
 \sum_{i=1}^k (A\t\tilde\x_i)\t (A\t\tilde\x_i).
 \label{update2}
 \end{equation}
 Furthermore
 \begin{equation}
 \lambda_i=(A\t\tilde\x_i)\t (A\t\tilde\x_i), i=1,\ldots,k.
 \label{lambdeq}
 \end{equation}
 \end{theo}

 We now explain the essence of Theorem \ref{updatel}.  View $\x_1,\ldots,\x_k$
 as approximation to the first $k$-left singular vectors of $A$, and
 $\sum_{i=1}^k \x_i (A\t \x_i)\t$ as the $k$-rank approximation to $A$.
 Hence $(A\t\x_i)\t (A\t \x_i)$ is an approximation to
 $\sigma_i^2$ of $A$ for $i=1,\ldots,k$.  Read additional vectors
 $\w_1,\ldots,\w_l\in\R^m$ such that at least one of this vectors
 is not in the subspace spanned by $\x_1,\ldots,\x_k$.
 Let $\X$ be the subspace spanned by $\x_1,\ldots,\x_k$ and
 $\w_1,\ldots,\w_l$. Hence $\x_1,\ldots,\x_{k},\ldots,\x_p$ is the orthonormal
 basis of $\X$ obtained from the vectors
 $\x_1,\ldots,\x_k, \w_1,\ldots,\w_l$, using the modified Gram-Schmidt
 algorithm.  (The modified Gram-Schmidt algorithm used to ensure the
 numerical stability.)
 Note that $k< p \le k+l$, and in general one
 has that $p=k+l$.  Consider the $p\times p$ nonnegative definite
 matrix $S=((A\t\x_i)\t (A\t \x_j))_{i,j=1}^p$.  Find its first
 $k$ eigenvectors to obtain $\tilde\x_1,\ldots,\tilde\x_k\in\R^m$ using (\ref{xtildef}).
 Then $C:=\sum_{i=1}^k \tilde\x_i (A\t\tilde \x_i)\t$ is the best
 approximation of $A$ by matrix $B$ of rank $k$ at most, whose columns are in the
 subspace $\X$.  In particular, the approximation $C$ is better
 than the previous approximation  $\sum_{i=1}^k \x_i (A\t \x_i)\t$,
 which is equivalent to the (\ref{update2}).

 $ \;$

 \noindent
 \textbf{Outline of Proof of Theorem \ref{updatel}.  }Let $S=(s_{ij})_{i,j=1}^p$.  Let
 $\e_i=(\delta_{i1},\ldots,\delta_{ip})\t,i=1,\ldots,p$ be the
 standard orthonormal basis in $\R^p$.  Use the
 definition of $A$ and Ky Fan characterization of the sum of the
 maximal $k$ eigenvalues of symmetric $S$ to deduce
 $$\sum_{i=1}^k(A\t\x_i)\t (A\t\x_i)=\sum_{i=1}^k \e_i\t S
 \e_i\le \sum_{i=1}^k \lambda_i=\sum_{i=1}^k \o_i\t S\o_i.$$
 Let $C:=A\t (\x_1,\ldots,\x_p)$.
 Then $S=C\t C$.  Hence the $\sqrt{\lambda_1}\ge\ldots\ge \sqrt{\lambda_p}$ are
 the singular values of $C$ and $\o_1,\ldots,\o_p$ are the right
 singular vectors of $C$.  Thus $C\o_i=\sqrt{\lambda_i}\bt_i\in \R^n$,
 where $\bt_i$ is the left singular vector of $C$ corresponding to
 the singular value $\sqrt{\lambda_i}$ for $i=1,\ldots,p$.
 A straightforward calculation shows that
 $\lambda_i=\o_i \t S\o_i = (A\t\tilde\x_i)\t (A\t\tilde\x_i)$ for  $i=1,\ldots,k$.
 Hence (\ref{update2}) holds.  Furthermore we also deduced the
 equalities in (\ref{lambdeq}) for $i=1,\ldots,k$.
 \qed

 \section{Algorithm}

 One starts the algorithm by choosing the first $k$-rank approximation
 to $A$ as follows.
 Let $\c_1,\ldots,\c_n\in \R^m$
 be the $n$ columns of $A$.
 Choose $k$ integers $1\le n_1<\ldots<n_k\le n$.
 Let $\x_1,\ldots,\x_q\in\R^m$
 be the orthonormal set obtained from $\c_{n_1},\ldots,\c_{n_k}$
 using the modified Gram-Schmidt algorithm.
 Set
 \begin{equation}\label{binit}
 B_0:=\sum_{i=1}^q \x_i(A\t\x_i)\t.
 \end{equation}
 In
 general $q=k$, but in some cases if $\x_1,\ldots,\x_k$ are
 linearly dependent, $q<k$.  Assume for simplicity of the
 exposition that $q=k$.  Then $B_0$ is of the form (\ref{bexpres})
 where $\y_i=A\t\x_i,i=1,\ldots,k$.

 Now update iteratively $k$-rank approximation of $B_{t-1}$ of $A$ to $B_t$, using
 Theorem \ref{updatel}, by letting
 $\w_1:=\c_{j_1},\ldots,\w_l:=\c_{j_l}\in\R^m$,
 for some $l$ integers
 $1\le j_1<\ldots<j_l\le n$.
 That is, we choose another $l$ sets of columns of $A$,
 preferably that were not chosen before, and update the
 $k$-rank approximation using the algorithm suggested by
 Theorem \ref{updatel} to obtain an improved $k$-rank
 approximation $B_t$ of $A$.

 Furthermore one can use the $k$-rank approximation $B_t$ from the above
 algorithm to approximate the first $k$-singular values
 $\sigma_1,\ldots,\sigma_k$, and
 the left and the right singular vectors
 $\u_1,\ldots,\u_k$ and $\v_1,\ldots,\v_k$ as follows.
 First, the square roots
 $\sqrt{\lambda_1(S)},\ldots,\sqrt{\lambda_k(S)}$ of the matrix
 $S$ are approximations for $\sigma_1,\ldots,\sigma_k$.
 Second,
 the vectors $\tilde\x_1,\ldots,\tilde\x_k$ approximate
 $\u_1,\ldots,\u_k$.  Let $\tilde \v_i: =A\t \tilde\x_i, i=1,\ldots,k$.
 Then $||\tilde \v_i||=\sqrt{\lambda_i(S)}$
 for $i=1,\ldots,k$.  Third, the renormalized $\tilde \v_i$ which are
 given as $\frac{1}{||\tilde \v_i||} \tilde \v_i$
 approximate the right singular eigenvectors
 $\v_i$ for $i=1,\ldots,k$.
 \[\]
 \framebox{\parbox[t]{5.4in}{

 \textbf{Fast $k$-rank approximation and SVD algorithm}

 \textbf{Input:}  positive integers $m,n,k,l,N$, $m\times n$ matrix
 $\bA$, $\epsilon > 0$.

 \noindent
 \textbf{Output:} an $m\times n$ $k$-rank approximation $B_f$ of
 $A$, with the ratios $\frac{||B_0||}{||B_t||}$ and $\frac{||B_{t-1}||}{||B_t||}$,
 approximations to $k$-singular values and $k$ left and right singular vectors of $A$.

 \noindent
 1.  Choose $k$-rank approximation $B_0$ using $k$ columns, (or rows), of $A$.

 \noindent
 2. \textbf{for} $t=1$ \textbf{to} $N$

 - Select $l$ columns, (or rows), from $A$ at random and update $B_{t-1}$ to $B_t$.

 - Compute the approximations to $k$-singular values, and $k$ left and right

 $\;\;$ singular vectors of $A$.

 - If $\frac{||B_{t-1}||}{||B_t||}> 1 - \epsilon$ let $f=t$ and
 finish.
 }}
 \[\]
 We now explain briefly the main steps of our algorithm.
 We read the dimensions $m,n$ of the data matrix $A$.  We set $N$
 as the maximal number of iterations we are going to execute to find
 the $k$-rank approximation of $A$.  We read the entries of the
 data matrix $A$, and finally the small parameter $\epsilon>0$.
 We choose the $k$-rank approximation $B_0$ using (\ref{binit}).
 Assume that $B_{t-1}$ is the current
 $k$-rank approximation to $A$.
 Next we choose additional $l$ columns of $A$ and
 update $B_{t-1}$ to $B_t$ using Theorem \ref{updatel}
 as explained in the previous section.  Recall that $||B_{t-1}||\le ||B_t||$.
 If the relative improvement in $||B_t||$ is less than $\epsilon$,
 i.e. $\frac{||B_{t-1}||}{||B_t||}> 1 - \epsilon$, we are
 satisfied with the approximation $B_t$ and finish our algorithm.
 If this does not happen then our algorithm stops after the $N$
 iteration.

 \section{Complexity of the algorithm}

 We assume here that $m\ge n$ and we consider our algorithm applied
 to the randomly selected columns of $A$.  (Otherwise we apply our
 algorithm to the rows of $A$.)
 The first step of our algorithm is to find an orthonormal
 $\x_1,\ldots,\x_q$ basis of the subspace spanned by
 $k$ randomly chosen columns of $A$.  The modified Gram-Schmidt
 algorithm needs $k^2 m$ flops assumed the worst, (generic), case $q=k$.
 This step needs $O(k^2m)$ flops
 \cite{GV}.
 Let $\y_i:=A\t\x_i\in\R^n$ for $i=1,\ldots,k$.  The computations
 of $\y_1,\ldots,\y_k$ needs $kmn$ operations and can be
 parallelized.  $B_0=\sum_{i=1}^k \x_i\y_i\t$ and does not have to
 be computed.

 To find $B_1$ we choose $l$ columns $\w_1,\ldots,\w_l$ of $A$ at
 random.  Then we
 perform the modified Gram-Schmidt algorithm on
 $\x_1,\ldots,\x_k,\w_1,\ldots,\w_l$ to find an orthonormal
 system $\x_1,\ldots,\x_p$, with $p\in [k,k+l]$.  This step requires
 $m(k+l)^2$ flops.  Let $\y_i=A\t\x_i$ for $i=k+1,\ldots,p$.
 This step needs at most $lmn$ flops and can be parallelized.
 Next we have the following two choices.
 First possibility is to find the SVD of $C_1=[\y_1,\ldots,\y_p]$.  This needs
 at most $O((k+l)^2 n)$ flops.
 Second possibility is to compute the matrix $S_1=C_1\t C_1$ given in
 Theorem \ref{updatel}.  This needs $nl(l+2k+1)/2$ flops.
 Then we compute the spectral decomposition of $S_1$ which needs
 $O((k+ l)^3)$ flops.  If $k+l$ is much smaller than $n$, it seems
 to us that the second method is more efficient.

 The first $k$ right singular vectors of $C_1$ are identical
 to the first $k$ eigenvectors of $S$.
 Use these $k$
 vectors, as explained in Theorem \ref{updatel}, to
 update $\x_1,\ldots,\x_k$ and $\y_1,\ldots,\y_k$. This needs $kp(m+n)$ flops.

 To update $B_{t-1}$ to $B_t$ for $t\ge 1$ require the same number
 of flops as to compute $B_1$.  Hence, the most intensive part of the
 computation lies in the computation of the spectral decomposition $S_t$,
 which is $O((k+l)^3)$.  Assuming that $l=O(k)$, we deduce that the intensive part
 of the computation is of order $O(k^3)$.  The simulations that
 we performed show that we need a small
 number of iterations, (around 5), to get a very good
 approximation to the best $k$-rank approximation of $A$.
 Hence our algorithm is expected to converge in $O(kmn)$ steps.

 \section{Simulation results}
 To assess the performance of our k-rank approximation algorithm we
 conducted different simulation on synthetic data and images. We
 also implemented our algorithm for three different sampling
 methods, uniform sampling with replacement, uniform sampling without
 replacement and weighted sampling for images according to weight of
 each row in gradient image. To show our algorithm guarantee
 convergence to the optimum (deterministic) k-rank approximation we
 applied the algorithm on a randomly generated data matrix of
 $3000\times 500$ with rank $500$. To measure the approximation
 error we defined the relative error of the approximation as,
 $\parallel A-B_t\parallel_{\rF}^2/\parallel A\parallel_{\rF}^2$,
 where $A $ is the original data matrix and $B_t$
 is k-rank approximation to $A$. Figure 1 shows convergence of the
 relative error to the optimum relative error ,
 $\|A-A_k\|_{\rF}^2/\|A\|_{\rF}^2$ where $A_k$ is the optimum true k-rank
 approximation to matrix $A$, as a function of iteration parameter
 $N $ for $k=100$ on the synthetic data matrix. In each iteration we
 randomly picked $l=100$ rows of the data matrix with and without
 replacement. For comparison the optimum relative error has been pointed
 out on the axis. As we expected the convergence property is faster when
 the rows are sampled without replacement. Figure 2 shows the same
 convergence result for the real images of Cameraman and Liftingbody (from Image Processing Toolbox of Matlab).
 Figure 3 shows
 the plots of the relative error versus total number of sampled
 rows in each method of sampling for the Cameraman and Liftingbody
 images.
 It can be seen that for these images our
 algorithm return very reasonable relative error only by using
 portion of the data.
 Needless to say that the importance of this iterative algorithm can
 be observed when it is applied to the big matrices of data where
 the SVD of the data matrix is often impossible to compute. In this
 case the simple randomized SVD methods \cite{FKV} may not be fast, since
 it needs to sample large enough number of rows to give the
 acceptable error bound on the approximation. To highlight this
 feature of our algorithm we applied our algorithm on a synthetic
 full rank data matrix of $8000\times 200$. We chose the parameter
 $l, N$ such that the relative error to be less than 2 times of optimum
 relative error. We observed that the speed up of our algorithm,
 which is the rate of the time needed by deterministic SVD to compute
 100-rank approximation to the time needed by our algorithm to compute
 $B_t$, is 42 in our system. The result for three more data sets with different values of $k$ has also
 been illustrated in Table 1. Due to system limitation and the comparison issue
 we were not able to show the speed up for bigger matrices where
 deterministic SVD (svd in Matlab) may fail to compute the svd of the
 matrix. However, our algorithm can always compute the k-rank approximation
 of the matrix as long as $k+l$ is not too big since the algorithm in each
 iteration deals with small matrices.

\begin{small}
\begin{table}
 \caption {Comparison of relative error and speed up of our algorithm with optimum $k$-rank approximation algorithm
 }
\begin{small}
 \begin{tabular}{|c|l|l|l|l|}
\hline

Data sets&Speed up & Re. ratio \\
 \hline
 Cameraman($256 \times 256$), $k=80$ &1.145& 1.083\\
\hline Liftingbody ($512 \times 512$), $k=100$&8 &1.08\\
\hline Map image($627 \times 865$) $k=200$&3.33 &1.067\\
\hline Random matrix($8000 \times 200$) $k=100$&42 &1.1\\

\hline
\end{tabular}
\end{small}

\end{table}
\end{small}

 \section{Conclusions}

 We have proposed a novel approach for fast computing of $k$-rank
 approximation of a given $m\times n$ data matrix, using
 Monte-Carlo method by choosing at random $l$ columns or rows of
 $A$.  The advantage of our algorithm is that we guarantee that
 every iteration improves the quality of our approximation. We
 applied our algorithm on synthetic data matrices as well as
 images and the result confirms the convergence of the
 relative error of the approximation to the optimum relative error.
 To highlight the important feature of this algorithm we applied this method
 on a big matrix of randomly generated data
 and we observed for the reasonable level of relative error the algorithm is
 also much faster than optimum k-rank approximation
 using deterministic SVD which may also fail to compute the svd for big matrices.
  We believe that this
 algorithm will have many applications in data mining, data storage and
 data analysis where dealing with high dimensional data is a major problem.

  \begin{figure}
 \centering
 \scalebox{0.5}{\includegraphics{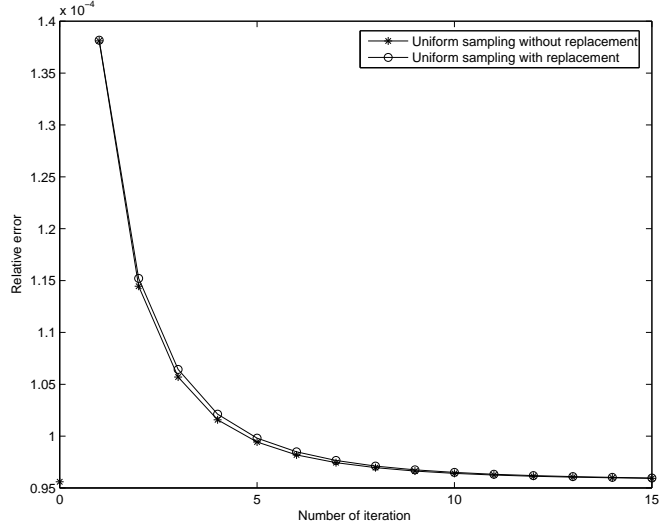}}
 \caption{Convergence property of the Monte-Carlo method for random data
 matrix($3000\times 500$).}\label{Fi:pic1}
 \end{figure}

 \begin{figure}
 \centering
 \mbox{\subfigure[Convergence property of the Monte-Carlo method for
 Cameraman image($256\times 256$), $k=80$.]
 {\epsfig{figure=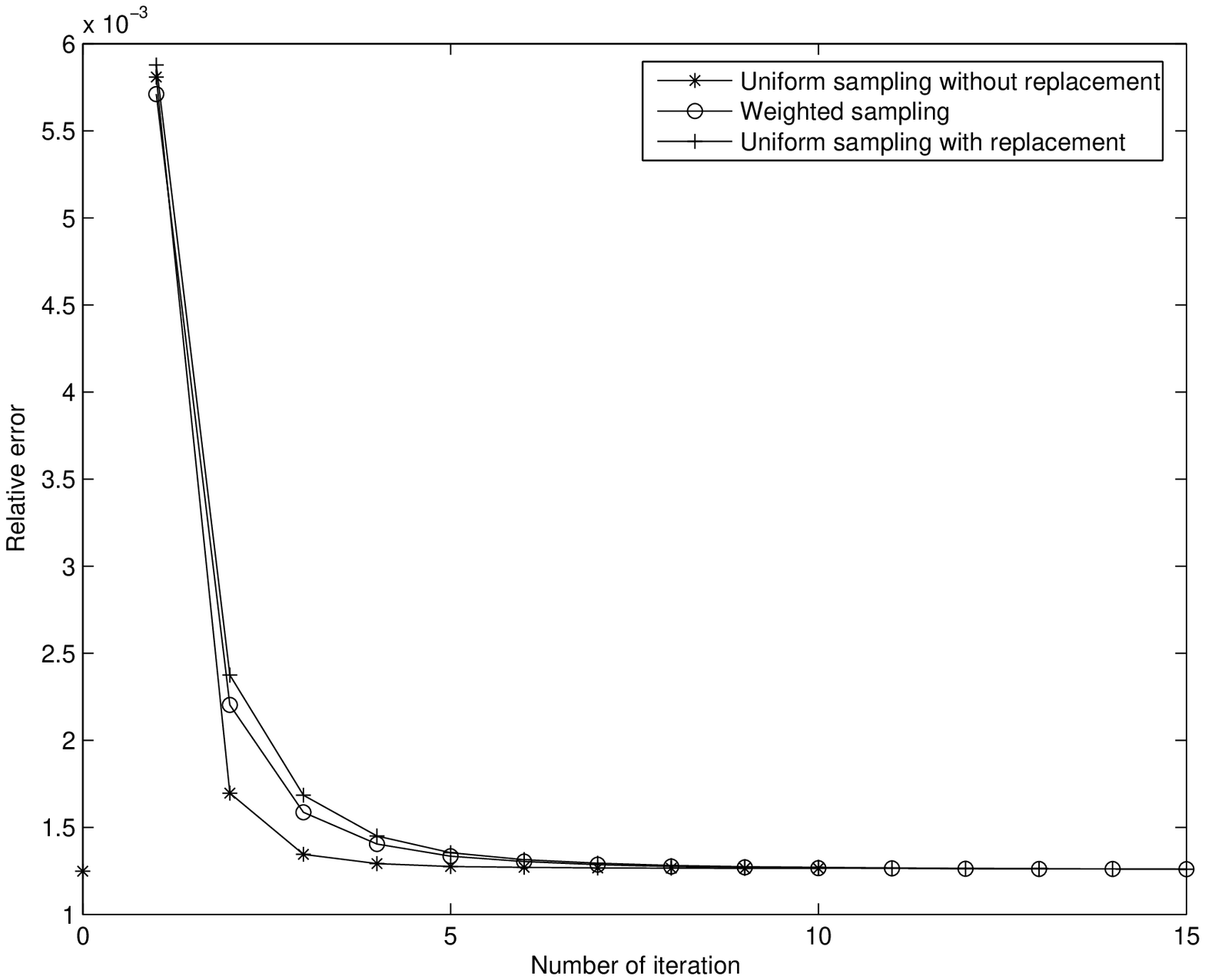,
 width=.50\textwidth}}\quad
 \subfigure[Convergence property of the Monte-Carlo method for
  Liftingbody image($512\times 512$), $k=80$.]
 {\epsfig{figure=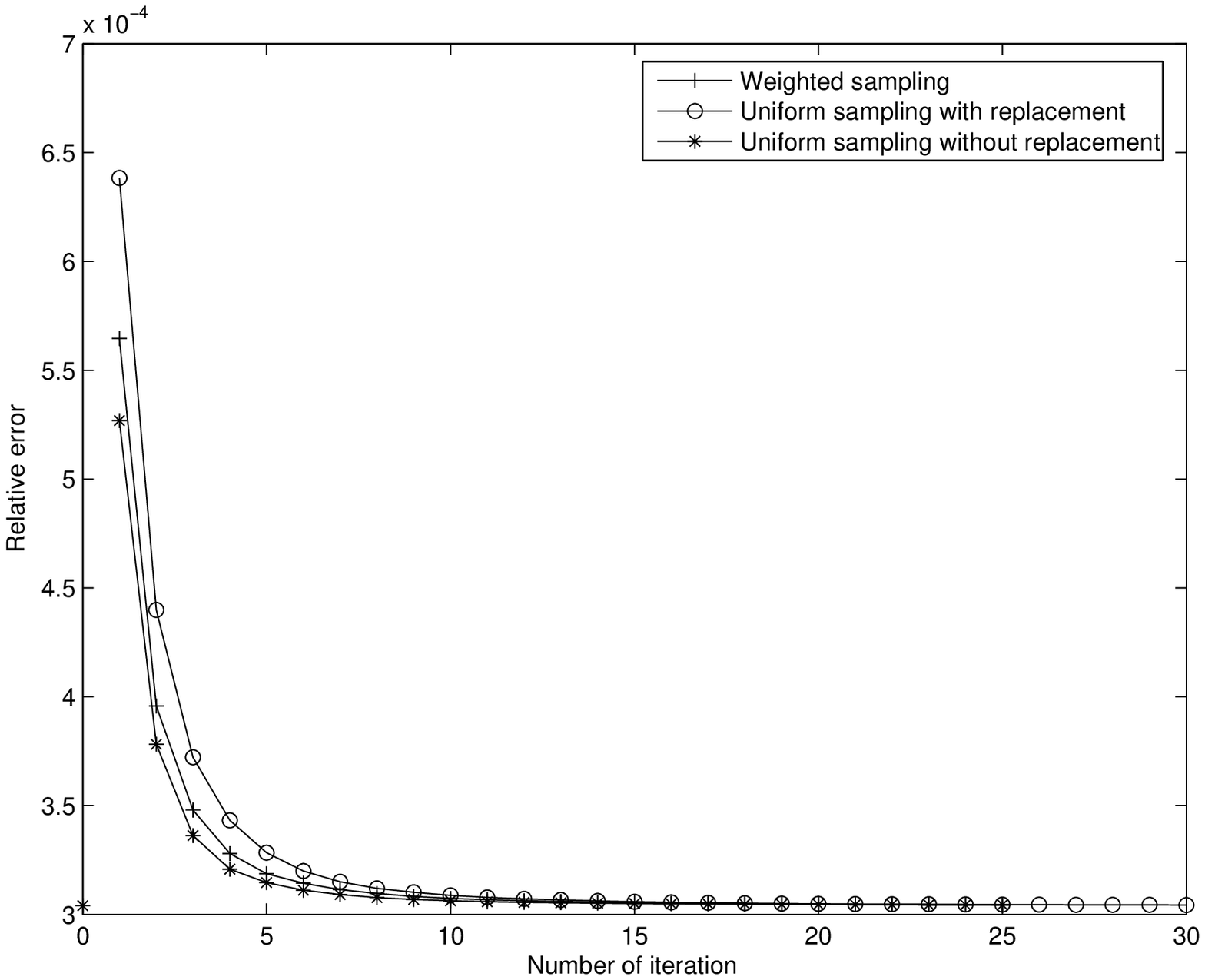,width=.50\textwidth}}}
 \caption{Cameraman and Liftingbody}\label{Fi:pic2}
 \end{figure}

 \begin{figure}
 \centering
 \mbox{\subfigure[Cameraman: Relative error versus total number
 of sampled rows, $k=80$.]
 {\epsfig{figure=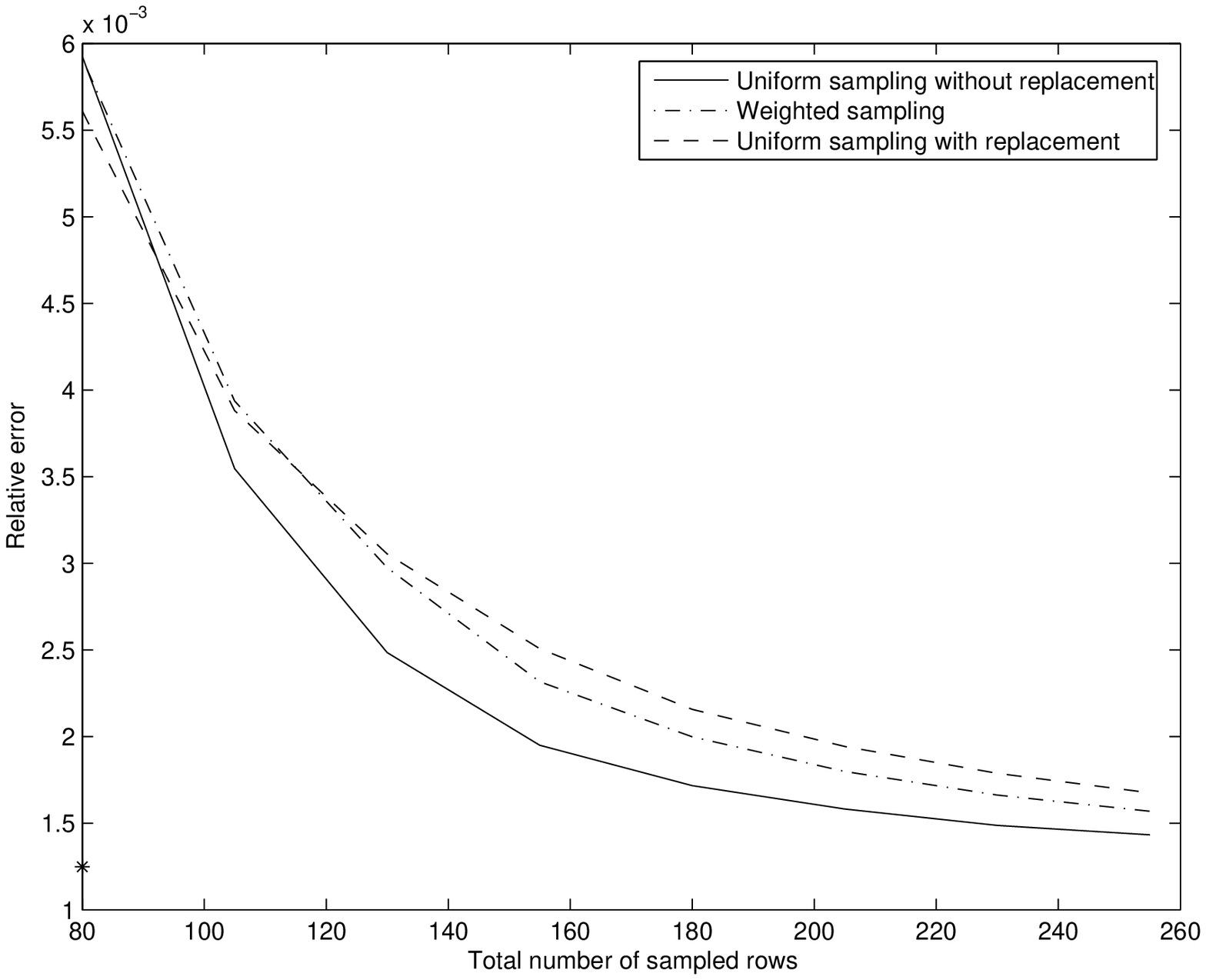,
 width=.50\textwidth}}\quad
 \subfigure[Liftingbody: Relative error versus total number of sampled
 rows, $k=80$.]
 {\epsfig{figure=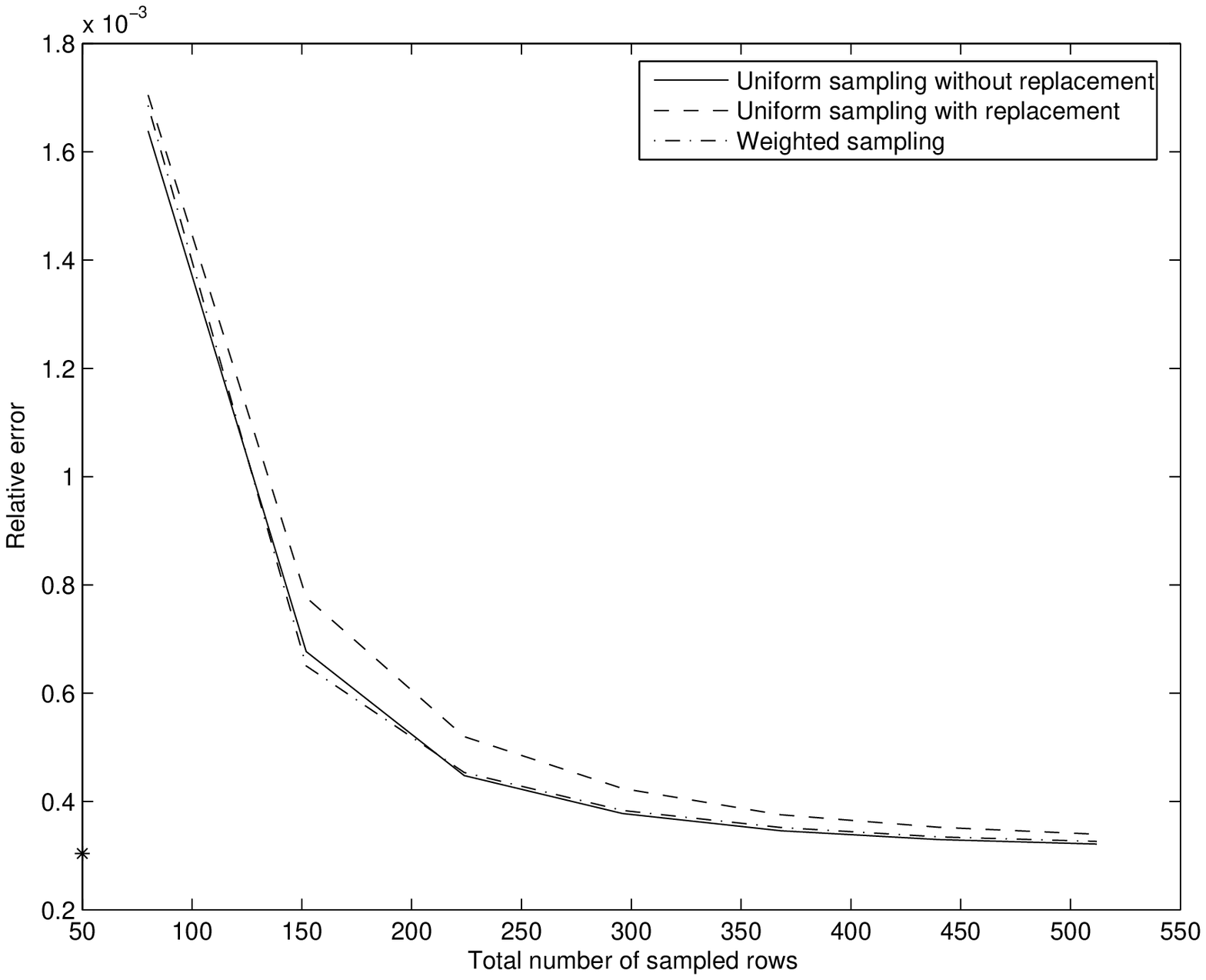,
 width=.50\textwidth}}}
 \caption{Cameraman and Liftingbody relative errors}\label{Fi:pic3}
 \end{figure}

\end{document}